# STATISTICAL PROPERTIES OF THE METHOD OF REGULARIZATION WITH PERIODIC GAUSSIAN REPRODUCING KERNEL

BY YI LIN[1] AND LAWRENCE D. BROWN[2]

*University of Wisconsin, Madison and University of Pennsylvania*

The method of regularization with the Gaussian reproducing kernel is popular in the machine learning literature and successful in many practical applications. In this paper we consider the periodic version of the Gaussian kernel regularization. We show in the white noise model setting, that in function spaces of very smooth functions, such as the infinite-order Sobolev space and the space of analytic functions, the method under consideration is asymptotically minimax; in finite-order Sobolev spaces, the method is rate optimal, and the efficiency in terms of constant when compared with the minimax estimator is reasonably high. The smoothing parameters in the periodic Gaussian regularization can be chosen adaptively without loss of asymptotic efficiency. The results derived in this paper give a partial explanation of the success of the Gaussian reproducing kernel in practice. Simulations are carried out to study the finite sample properties of the periodic Gaussian regularization.

**1. Introduction.** The method of regularization is a popular approach for nonparametric function estimation. Let $f$ be the nonparametric function to be estimated. The method of regularization takes the form

$$(1) \qquad \min_{f \in \mathcal{F}}[L(f, \text{data}) + \lambda J(f)],$$

where $L$ is the empirical loss, often taken to be the negative log-likelihood, and $J(f)$ is the penalty functional, usually a quadratic functional corresponding to a norm or semi-norm of a reproducing kernel Hilbert space $\mathcal{F}$. Most often the penalty functional is chosen so that smoother functions

Received August 2002; revised July 2003.
[1]Supported in part by NSF Grant DMS-01-34987.
[2]Supported in part by NSF Grant DMS-99-71751.

incur smaller penalty. The smoothing parameter $\lambda$ controls the tradeoff between minimizing the empirical loss and obtaining a smooth solution. For a concrete example, let us look at the regression model

$$(2) \qquad y_j = f(x_j) + \delta_j, \qquad j = 1, \ldots, n,$$

where $x_j \in R$, $j = 1, \ldots, n$, are the regression inputs, $y_j$'s are the responses, and $\delta_j$'s are independent $N(0,1)$ noises. In this case we may take $L(f, \text{data}) = \sum_{j=1}^{n} (y_j - f(x_j))^2$ in the method of regularization (1).

The reproducing kernel Hilbert space $\mathcal{F}$ is typically of infinite dimension. In many situations, including regression and generalized regression, when the penalty functional $J(f)$ is a norm over $\mathcal{F}$, the representer theorem [Kimeldorf and Wahba (1971)] guarantees that the solution to (1) over $\mathcal{F}$ falls in the finite-dimensional space spanned by $\{K(x_j, \cdot), j = 1, \ldots, n\}$, where $K(\cdot, \cdot)$ is the reproducing kernel corresponding to $J(f)$. See also Schölkopf, Herbrich and Smola (2001) for some generalizations of the representer theorem. Therefore, we can write the solution as $\hat{f} = \sum_{j=1}^{n} c_j K(x_i, x)$. The minimization problem can then be solved in this finite-dimensional space.

The smoothing spline well known in the nonparametric statistics literature is an example of the method of regularization. In the smoothing spline the reproducing kernel Hilbert space $\mathcal{F}$ is a Hilbert Sobolev space and the penalty functional $J(f)$ is the norm or semi-norm of the space, such as $\int [f^{(m)}(x)]^2 \, dx$. The commonly used cubic smoothing spline corresponds to the case $m = 2$. The reproducing kernel of the Hilbert Sobolev space was given in Wahba (1990).

The method of regularization has also been popular in the machine learning literature. Examples include regularization networks and more recently, support vector machines. See, for example, Girosi, Jones and Poggio (1993), Smola, Schölkopf and Müller (1998), Wahba (1999) and Evgeniou, Pontil and Poggio (2000). One reproducing kernel that is particularly popular in the machine learning literature is the Gaussian reproducing kernel (commonly referred to as the Gaussian kernel in the machine learning literature, not to be confused with the Gaussian kernel used in kernel smoothing in the nonparametric statistics literature). Let $G(r) = (2\pi)^{-1/2} \omega^{-1} \exp(-r^2/(2\omega^2))$ be the density function of $N(0, \omega^2)$. The Gaussian reproducing kernel has the form $G(s, t) \equiv G(s - t)$. This is a common example of the translation invariant reproducing kernels popular in machine learning. It is known [Girosi, Jones and Poggio (1993) and Smola, Schölkopf and Müller (1998)] that the Gaussian reproducing kernel corresponds to the penalty functional (up to a constant)

$$(3) \qquad J_g(f) = \sum_{m=0}^{\infty} \frac{\omega^{2m}}{2^m m!} \int_{-\infty}^{\infty} [f^{(m)}(x)]^2 \, dx.$$

Smola, Schölkopf and Müller (1998) introduced the periodic Gaussian reproducing kernel for estimating $2\pi$-periodic functions in $[-\pi, \pi]$ as the reproducing kernel corresponding to the penalty functional

$$(4) \qquad J_0(f) = \sum_{m=0}^{\infty} \frac{\omega^{2m}}{2^m m!} \int_{-\pi}^{\pi} [f^{(m)}(x)]^2 \, dx.$$

From (3) and (4) it is clear that the two reproducing kernels are closely related. The connection between the two reproducing kernels will be clearer when we consider the computation with the periodic Gaussian reproducing kernel in Section 5.

Many researchers in machine learning have derived upper bounds of the generalization performance of the method of regularization with the Gaussian or periodic Gaussian reproducing kernels. See Williamson, Smola and Schölkopf (2001) and the references therein. However, while popular in the machine learning literature, and successful in many practical applications, the statistical asymptotic properties of the method of regularization with the Gaussian or periodic Gaussian reproducing kernels have not been studied systematically. In this paper we study the asymptotic properties of the method of regularization with the periodic Gaussian reproducing kernel in nonparametric function estimation problems and derive the asymptotic risk (up to constants) of the method of regularization with the periodic Gaussian reproducing kernel. We choose to work with the periodic Gaussian reproducing kernel because it allows a detailed asymptotic analysis. We believe the results obtained in this paper should also give insights on the statistical properties of the Gaussian reproducing kernel.

Motivated by the equivalence results of Brown and Low (1996) for Gaussian nonparametric regression and Nussbaum (1996) for density estimation [see also Golubev and Nussbaum (1998) for spectral density estimation; Grama and Nussbaum (1997) for nonparametric generalized linear regression], we first look at the white noise problem

$$(5) \qquad Y_n(t) = \int_{-\pi}^{t} f(u) \, du + n^{-1/2} B(t), \qquad t \in [-\pi, \pi],$$

where $B(t)$ is a standard Brownian motion on $[-\pi, \pi]$ and we observe $Y_n = (Y_n(t), -\pi \leq t \leq \pi)$. We consider the situation where the function $f$ belongs to a certain function ellipsoid of the form

$$(6) \qquad \left\{ f : f(t) = \sum_{l=0}^{\infty} \theta_l \phi_l(t), \sum_{l=0}^{\infty} \rho_l \theta_l^2 \leq Q \right\},$$

for some positive sequence $\{\rho_l, l = 0, 1, \ldots\}$. Here $\{\phi_0(t) = (2\pi)^{-1/2}, \phi_{2l-1}(t) = \pi^{-1/2} \sin(lt), \phi_{2l}(t) = \pi^{-1/2} \cos(lt)\}$ is the classical trigonometric basis in $L_2(-\pi, \pi)$ and $\theta_l = (f, \phi_l)$ is the corresponding Fourier coefficient, where $(f, \phi) = \int_{-\pi}^{\pi} f(t)\phi(t)\,dt$ denotes the usual inner product in $L_2(-\pi, \pi)$.

The commonly considered Sobolev ellipsoid $H^m(Q)$ corresponds to the sequence $\rho_0 = 1$, $\rho_{2l-1} = \rho_{2l} = l^{2m} + 1$ in (6). This is the $m$th order Sobolev space of periodic functions on $[-\pi, \pi]$. An alternative definition of $H^m(Q)$ is

$$
(7) \quad H^m(Q) = \bigg\{ f \in L^2(-\pi, \pi) : f \text{ is } 2\pi\text{-periodic}, \\
\int_{-\pi}^{\pi} [f(t)]^2 + [f^{(m)}(t)]^2 \, dt \leq Q \bigg\}.
$$

Therefore, the $m$th order Sobolev space consists of functions that possess $m$th order smoothness. The order of smoothness is determined by the rate at which the sequence of $\rho$'s increases. In the Sobolev space case the rate is of polynomial order.

Another function space that has been considered in the literature is the space of analytic functions. An ellipsoid of analytic functions $A_\alpha(Q)$ corresponds to (6) with the exponentially increasing sequence $\rho_l = \exp(\alpha l)$, where $\alpha$ is a positive constant. Such a function space can be motivated by considering the Fourier series in complex exponentials and considering the domain in which the function is analytical. For details, see Johnstone (1998). The norm of this function space can not be expressed in terms of integrals of squared derivatives of integer order.

We now introduce a new function space $H_\omega^\infty$ that can be seen as the Sobolev space of infinite order,

$$
(8) \quad H_\omega^\infty(Q) = \bigg\{ f : f(t) = \sum_{l=0}^{\infty} \theta_l \phi_l(t), \sum_{l=0}^{\infty} \rho_l \theta_l^2 \leq Q; \\
\rho_0 = 1, \rho_{2l-1} = \rho_{2l} = e^{l^2 \omega^2/2} \bigg\},
$$

where $\omega$ is a positive constant, and $\phi$'s are the classical trigonometric basis over $(-\pi, \pi)$. Simple calculation shows that an equivalent definition of $H_\omega^\infty(Q)$ is

$$
H_\omega^\infty(Q) = \bigg\{ f \in L^2(-\pi, \pi) : f \text{ is } 2\pi\text{-periodic}, \\
\sum_{m=0}^{\infty} \frac{\omega^{2m}}{m! 2^m} \int_{-\pi}^{\pi} [f^{(m)}(x)]^2 \, dx \leq Q \bigg\}.
$$

From this we can see that $H_\omega^\infty$ can be seen as the Sobolev space of infinite order, and that the penalty functional $J_0$ of the periodic Gaussian reproducing kernel as defined in (4) corresponds to the norm of $H_\omega^\infty(Q)$.

In this paper we focus on the method of regularization with the periodic Gaussian penalty (4). We will refer to this method as periodic Gaussian regularization. We study the statistical properties of this method both in the situation that $f \in H_\omega^\infty$ and the situation $f \notin H_\omega^\infty$.

By converting the functions into the corresponding sequence of Fourier coefficients, we can see that the white noise problem (5) is equivalent to the following Gaussian sequence model:

$$(9) \qquad y_l = \theta_l + \varepsilon_l, \qquad l = 0, 1, \ldots,$$

where the $\varepsilon_l$'s are independent $N(0, 1/n)$ noises and the $\theta_l$'s are the Fourier coefficients of $f$. The periodic Gaussian regularization corresponds to

$$(10) \qquad \min \sum_{l=0}^\infty (y_l - \theta_l)^2 + \lambda \sum_{l=0}^\infty \beta_l \theta_l^2$$

with $\beta_l = e^{l^2 \omega^2/2}$.

In Section 2 we establish the asymptotic minimax risk (up to the constant) of nonparametric problems in the space $H_\omega^\infty(Q)$, and show that the periodic Gaussian regularization achieves this optimal asymptotic risk. In Section 3 we study the asymptotic performance of the periodic Gaussian regularization in the situation where the underlying function to be estimated is in the Sobolev ellipsoid $H^m(Q)$ with unknown $m$ and $Q$, or in the analytic function ellipsoid $A_\alpha(Q)$ with unknown $\alpha$ and $Q$. We show that the method under study is asymptotically minimax in analytic function ellipsoids. For Sobolev ellipsoids $H^m(Q)$, the periodic Gaussian regularization achieves the optimal rate of convergence, and the efficiency in terms of the constant is reasonably high, tending to 1 as $m$ goes to infinity.

In Section 4 we consider choosing the smoothing parameters with the unbiased estimator of risk. The procedure is the well known Mallows' $C_p$ [Mallows (1973)], sometimes called Mallows' $C_L$ in the literature. Li (1986, 1987) established the asymptotic optimality of $C_p$ in many nonparametric function estimation methods, including the method of regularization. Kneip (1994) obtained oracle inequalities for choosing smoothing parameters with $C_p$ in ordered linear smoothers. See also Cavalier, Golubev, Picard and Tsybakov (2002). These results can be used to study the periodic Gaussian regularization with smoothing parameters chosen by the unbiased risk estimator. We show that the resulting data-driven method retains the good theoretical properties of the periodic Gaussian regularization established in Sections 2 and 3. Thus, adaptive estimation is achieved for unknown order of smoothness by the periodic Gaussian regularization in the white noise model.

Due to the equivalence between the white noise model and other statistical models, we expect the periodic Gaussian regularization to have good statistical properties in other situations such as regression and generalized regression. In fact, the equivalence results in Brown and Low (1996) show that the asymptotic results we obtained in Sections 2–4 for the white noise model apply to the periodic Gaussian regularization in the regression problem (2) with fixed equidistant design. In regression problems with nonequidistant design, the periodic Gaussian regularization in regression does not match up exactly with the periodic Gaussian regularization in the white noise model, and therefore our results do not translate directly. However, we believe the results in the white noise model still give insights to the regression problem with general design. In this connection, see Brown and Zhao (2002).

In Section 5 we consider the computation of the periodic Gaussian regularization in regression. The computation does not require equidistant design. Some simulations are given in Section 6 to study the finite sample properties of the periodic Gaussian regularization. In particular, the effect of the joint tuning of the smoothing parameters is studied, and the periodic Gaussian regularization is compared with the periodic cubic smoothing spline on four functions of different orders of smoothness. The simulation suggests that the finite sample performance of the periodic Gaussian regularization is comparable to that of the periodic cubic smoothing spline when the regression function is of moderate smoothness. In the case of a very smooth function, the periodic Gaussian regularization may have an advantage. Summary and discussion are given in Section 7. Technical proofs are relegated to Section 8.

Throughout this paper the expression $a_n \sim b_n$ means that $a_n/b_n \to 1$ as $n \to \infty$.

**2. Estimation in the Sobolev space of infinite order.** In this section we consider the white noise problem in $H_\omega^\infty(Q)$.

THEOREM 1. *The asymptotic minimax risk for nonparametric function estimation in the infinite-order Sobolev ellipsoid $H_\omega^\infty(Q)$ is $2\sqrt{2}\omega^{-1}n^{-1}(\log n)^{1/2}$. That is,*

$$\inf_{\bar\theta} \sup_{\theta \in H_\omega^\infty(Q)} \sum_{l=0}^\infty E(\bar\theta_i - \theta_i)^2 \sim 2\sqrt{2}\omega^{-1}n^{-1}(\log n)^{1/2},$$

*where the infimum is over all possible estimators $\bar\theta$.*

Notice this asymptotic minimax risk does not depend on $Q$, but depends on $\omega$.

In the following we consider the periodic Gaussian regularization. The following lemma will be used several times in later proofs.

LEMMA 1. *Consider the periodic Gaussian regularization* (10) *in the white noise model. Denote the estimator by* $\hat{\theta}$. *We have* $\sum \mathrm{var}\,\hat{\theta} \sim 2\sqrt{2}\omega^{-1}n^{-1} \times (-\log \lambda)^{1/2}$, *as* $n \to \infty$ *and* $\lambda_{(n)} \to 0$.

THEOREM 2. *The periodic Gaussian regularization* (10) *in the white noise model is asymptotically minimax in the infinite-order Sobolev ellipsoid* $H_\omega^\infty(Q)$, *if the smoothing parameter* $\lambda$ *satisfies*

(11) $$\log(1/\lambda) \sim \log n \quad \text{and} \quad \lambda = o(n^{-1}(\log n)^{1/2}).$$

*That is,*

$$\inf_\lambda \sup_{\theta \in H_\omega^\infty(Q)} \sum_{l=0}^\infty E(\hat{\theta}_i - \theta_i)^2 \sim 2\sqrt{2}\omega^{-1}n^{-1}(\log n)^{1/2},$$

*and this asymptotic risk is achieved when* (11) *is satisfied. Here* $\hat{\theta}$ *is the method of regularization estimator from* (10) *with* $\beta_l = e^{l^2\omega^2/2}$.

The condition (11) is satisfied if $n\lambda_n$ is bounded away from zero and infinity, but is milder. For example, it is satisfied by sequences $\lambda_n = Cn^{-1}(\log n)^\alpha$ for any constants $C > 0$ and $-\infty < \alpha < 1/2$. The adaptive choice of $\lambda$ is considered in Section 4.

**3. Estimation over Sobolev spaces and spaces of analytic functions.** In this section we consider the performance of the periodic Gaussian regularization when the function $f$ to be estimated in the white noise problem belongs to a Sobolev body $H^m(Q)$ with unknown $m$ and $Q$, or an analytic function ellipsoid $A_\alpha(Q)$ with unknown $\alpha$ and $Q$. In these cases the function to be estimated does not lie in the function space used in the method of regularization.

THEOREM 3. *Assume* $f \in H^m(Q)$ *with* $m \geq 1$ *in the white noise model* (5). *Consider the periodic Gaussian regularization estimator* $\hat{\theta}$ (10) *with* $\beta_{2l-1} = \beta_{2l} = \exp(l^2\omega^2/2)$. *We have*

$$\inf_\lambda \sup_{\theta \in H^m(Q)} \sum_l E(\hat{\theta}_l - \theta_l)^2 \sim (2m+1)m^{-2m/(2m+1)}Q^{1/(2m+1)}n^{-2m/(2m+1)}.$$

*This asymptotic risk is achieved when* $\log(1/\lambda)/\omega^2 \sim (mnQ)^{2/(2m+1)}/2$.

REMARK 1. *The conclusion of Theorem 3 holds for noninteger* $m > 1$.

For the ellipsoid $A_\alpha(Q)$ of analytic functions, we have the following:

THEOREM 4. *Assume $f \in A_\alpha(Q)$ in the white noise problem* (5). *Consider the periodic Gaussian regularization estimator $\hat{\theta}$ from* (10) *with $\beta_{2l-1} = \beta_{2l} = e^{l^2 \omega^2/2}$. We have*

$$\inf_\lambda \sup_{\theta \in A_\alpha(Q)} \sum_l E(\hat{\theta}_l - \theta_l)^2 \sim 2n^{-1}\alpha^{-1} \log n.$$

*This asymptotic risk is achieved when $\log(1/\lambda)/\omega^2 = (\log n)^2/(2\alpha^2)$.*

The proof of this theorem is similar to that of Theorem 3, with $\rho_l = e^{\alpha l}$, and is skipped. It is known that the asymptotic minimax risk in $A_\alpha(Q)$ is $2n^{-1}\alpha^{-1} \log n$; see Johnstone (1998). Therefore, Theorem 4 says that the periodic Gaussian regularization is asymptotically minimax in $A_\alpha(Q)$.

We can study the asymptotic efficiency of the periodic Gaussian regularization compared with the minimax estimator for nonparametric problems in $H^m(Q)$. We consider the maximum asymptotic risk over $H^m(Q)$. We compare the minimum of such asymptotic risk achieved by the periodic Gaussian regularization with the minimax risk over $H^m(Q)$. This indicates how close to the minimax value one can get with the periodic Gaussian regularization. A similar study had been carried out by Carter, Eagleson and Silverman (1992), who studied the efficiency of the cubic smoothing spline in the second-order Sobolev space.

It is well known that the asymptotic minimax risk over $H^m(Q)$ is

$$[2m/(m+1)]^{2m/(2m+1)}(2m+1)^{1/(2m+1)}Q^{1/(2m+1)}n^{-2m/(2m+1)}.$$

This can be derived with an argument along the line of the proof of Theorem 1. Figure 1, left panel, gives the ratio between the asymptotic risk of the periodic Gaussian regularization and the minimax risk when the sample size $n$ is kept to be the same. The right panel gives the efficiency of the periodic Gaussian regularization. The efficiency is calculated in terms of sample sizes needed to achieve the same risk. We can see that the efficiency goes to one when the function is very smooth. The lowest efficiency occurs when $m = 1$, and the lowest efficiency is 33.3%. The efficiency when $m = 2$ is 53.3%.

**4. Adaptive choice of the smoothing parameter.** In the earlier sections we studied the performance of the periodic Gaussian regularization when the smoothing parameter $\lambda$ has an appropriate rate of decrease. This appropriate rate depends on $m$ (or $\alpha$ or $\omega$) and $Q$, which are generally unknown in practice. In this section we consider the problem of choosing the smoothing parameter with data. We study the common approach of choosing the smoothing parameter through the unbiased estimator of risk (Mallows' $C_p$). By making use of the oracle inequalities developed in Kneip (1994) [see also

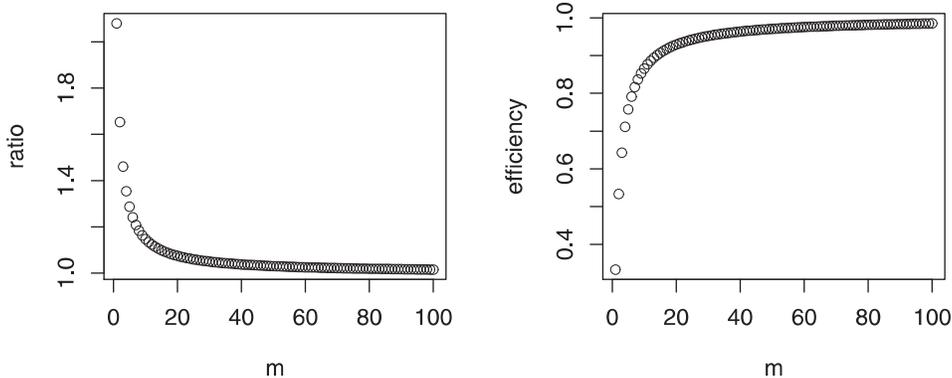

Fig. 1. *The efficiency of the periodic Gaussian regularization method.*

Cavalier, Golubev, Picard and Tsybakov (2002)], we show that the estimator chosen by the unbiased estimator of risk has the same asymptotic risk as the estimator with the optimal (theoretical) smoothing parameter. Thus, no asymptotic efficiency is lost due to not knowing $m$, $Q$ and $\omega$.

The number $\omega$ appears in the asymptotic risk of the periodic Gaussian regularization estimator in the function space $H_\omega^\infty(Q)$, but does not play an important role in the asymptotic risk in the function space $H^m(Q)$, so long as $\lambda$ is suitably chosen. From (22) in the proof of Theorem 3 we can see that the leading terms in the asymptotic risk in $H^m(Q)$ depend on $\omega$ and $\lambda$ only through $-\log\lambda/\omega^2$. The asymptotic results suggest that tuning one of $\lambda$ and $\omega$ may suffice. For finite sample size, though, it may pay to tune $\omega$ as well as $\lambda$. Usually there is a range of $\omega$ that works almost equally well if $\lambda$ is tuned correspondingly and vice versa. See the simulation in Section 6 for examples. Thus, we consider a rough tuning for $\omega$, just to get to a reasonable range, and a fine tuning over $\lambda$.

Formally, we take a finite number of $\omega$'s: $\omega_1,\ldots,\omega_S$, and tune $\lambda$ and $\omega$ jointly over $\lambda$ and $\omega_s \in \{\omega_1,\ldots,\omega_S\}$. For asymptotic consideration, a range of $[0,1]$ for $\lambda$ suffices, since asymptotically $\lambda$ should go to zero. In practice we may use a slightly larger range.

The tuning is based on the unbiased estimator of risk. Writing

$$\tau_l = (1 + \lambda\beta_l)^{-1},$$

our estimator is

$$\hat{\theta}_l = \tau_l y_l.$$

We can express the risk of our estimator as

$$\sum_l E(\hat{\theta}_l - \theta_l)^2 = (1/n)\sum_{l=0}^\infty \tau_l^2 + \sum_{l=0}^\infty (1-\tau_l)^2 \theta_l^2.$$

Now an unbiased estimator for $\theta_l^2$ is $y_l^2 - (1/n)$. Plugging in, we get that

$$(12) \quad \sum_{l=0}^{\infty}[(\tau_l^2 - 2\tau_l)(y_l^2 - 1/n) + (1/n)\tau_l^2] = \sum_{l=0}^{\infty}[(\tau_l^2 - 2\tau_l)y_l^2 + (2/n)\tau_l]$$

is an unbiased estimator of $\sum_l E(\hat{\theta}_l - \theta_l)^2 - \sum \theta_l^2$. We choose $\lambda^*$ and $\omega^*$ that minimize the unbiased risk (12), and use the corresponding periodic Gaussian regularization estimator $\hat{\theta}^*$. Kneip (1994) studied the adaptive choice among ordered linear smoothers with the unbiased risk estimator. A family of ordered linear smoothers satisfies the condition that for any member $\hat{\theta}_l = \tau_l y_l$, $l = 0, 1, \ldots$, of the family, we have $\tau_l \in [0, 1]$ $\forall l$; and for any two members of the family, $\tau_l y_l$ and $\tau_l' y_l$, $l = 0, 1, \ldots$, we have either $\tau_l \geq \tau_l'$ $\forall l$, or $\tau_l' \geq \tau_l$ $\forall l$. It is easy to check that for any fixed $\omega \in \{\omega_1, \ldots, \omega_S\}$, the method of regularization estimators with varying $\lambda$ form a family of ordered linear smoothers. Applying the result in Kneip (1994) [recast in the Gaussian sequence model setting in Cavalier, Golubev, Picard and Tsybakov (2002)] to our situation gives the following:

LEMMA 2. *Consider the Gaussian sequence model* (9) *and the periodic Gaussian regularization* (10). *Suppose $\lambda^*$ and $\omega^*$ minimize* (12) *over $\lambda \in [0, 1]$ and $\omega \in \{\omega_1, \ldots, \omega_S\}$, and $\hat{\theta}^*$ is the corresponding periodic Gaussian regularization estimator. Then there exist positive constants $C_1$ and $C_2$ such that for any $\theta \in l^2$ and any positive constant $B$, we have*

$$(13) \quad \sum_l E(\hat{\theta}_l^* - \theta_l)^2 \leq (1 + C_1 B^{-1}) \min_{\lambda, \omega_s} \left\{ \sum_l E(\hat{\theta}_l - \theta_l)^2 \right\} + n^{-1} C_2 B.$$

We then have the following:

THEOREM 5. *For the periodic Gaussian regularization estimator $\hat{\theta}^*$ chosen by the unbiased estimator of risk, we have*

$$\sup_{\theta \in H_{\omega_s}^{\infty}(Q)} \sum_l E(\hat{\theta}_l^* - \theta_l)^2 \sim 2\sqrt{2}\omega_s^{-1} n^{-1}(\log n)^{1/2} \qquad \forall s \in \{1, \ldots, S\},$$

$$\sup_{\theta \in H^m(Q)} \sum_l E(\hat{\theta}_l^* - \theta_l)^2 \sim (2m+1)m^{-2m/(2m+1)} Q^{1/(2m+1)} n^{-2m/(2m+1)},$$

$$\sup_{\theta \in A_\alpha(Q)} \sum_l E(\hat{\theta}_l^* - \theta_l)^2 \sim 2n^{-1}\alpha^{-1} \log n.$$

Therefore, the adaptive periodic Gaussian regularization estimator $\hat{\theta}^*$ is asymptotically minimax in $H_{\omega_s}^{\infty}(Q)$ and $A_\alpha(Q)$, and achieves the optimal rate in $H^m(Q)$. The asymptotic efficiency is the same as that given in Section 3. Hence, the estimator adapts to any unknown order of smoothness.

**5. Computation of periodic Gaussian regularization in regression.** In order for the periodic Gaussian regularization in regression and generalized regression to be practically computable, we need the form of the reproducing kernel corresponding to the penalty functional $J_0(f)$, that is, the reproducing kernel of $H_\omega^\infty$. Smola, Schölkopf and Müller (1998) gave the following expression for the periodic Gaussian reproducing kernel:

$$(14) \qquad R(s,t) = (1/\pi) \sum_{l=1}^{\infty} \exp(-l^2\omega^2/2)\cos(l(s-t)).$$

Due to the fast decay of the sequence $\exp(-l^2\omega^2/2)$, it is possible to approximate the series (14) with finitely many terms. However, an alternative formula of the kernel (14) is better suited for computation. We first state a lemma due to Williamson, Smola and Schölkopf (2001).

LEMMA 3. *Let $V(s-t)$ be a reproducing kernel with $V: R \to R$ being an even function. Let*

$$V_\nu(s) = \sum_{k=-\infty}^{\infty} V(s-k\mu).$$

*Then*

$$V_\nu(s-t) = \frac{\sqrt{2\pi}}{\nu}\tilde{V}(0) + \sum_{k=1}^{\infty} \frac{2}{\nu}\sqrt{2\pi}\tilde{V}\left(\frac{2k\pi}{\nu}\right)\cos\frac{2k\pi(s-t)}{\nu},$$

*where $\tilde{V}$ is the Fourier transform of $V$.*

Define $G^\infty(r) = \sum_{k=-\infty}^{\infty} G(r-2k\pi)$. It follows directly from Lemma 3 that $G^\infty(s-t)$ is the reproducing kernel (14) corresponding to the periodic Gaussian regularization. The function $G^\infty$ can be approximated with the finite series $G^J = \sum_{k=-J}^{J} G(s-2k\pi)$ for some $J$. In fact, we have

$$0 < G^\infty(s) - G^1(s) < 2.1 \times 10^{-20} \qquad \forall s \in [-\pi, \pi] \text{ for } \omega \le 1.$$

For $\omega > 1$, we can choose a positive integer $J$ such that $2J+1 \ge 3\omega$. Then $0 < G^\infty(s) - G^J(s) < 10^{-20}\ \forall s \in [-\pi, \pi]$. Therefore, $G^J(s)$ is an easily computable proxy of $G^\infty(s)$.

Now consider the periodic Gaussian regularization in the regression problem (2) with the empirical loss being $\sum_{j=1}^{n}(y_j - f(x_j))^2$. Here we assume $x_j \in (-\pi, \pi)$, $j = 1, \ldots, n$, and the regression function $f$ is $2\pi$-periodic. The theory of reproducing kernel Hilbert spaces guarantees that the solution to the method of regularization falls in a finite-dimensional space spanned by $G^\infty(x_j, \cdot)$. That is, we can write $\hat{f}(x) = \sum_{j=1}^{n} \hat{c}_j G^\infty(x_j - x)$, and the penalized regression (1) becomes

$$(y - G^\infty c)'(y - G^\infty c) + \lambda c' G^\infty c,$$

where, with little risk of confusion, we write $y = (y_1, \ldots, y_n)^t$, $c = (c_1, \ldots, c_n)^t$, and $G^\infty$ is the $n \times n$ matrix $(G^\infty(x_i - x_j))$. The solution can then be found to be $\hat{c} = (G^\infty + \lambda I)^{-1} y$. In order to compute the solution as well as Mallows' $C_p$ for tuning the smoothing parameters, we use the eigenvalue–eigenvector decomposition $G^\infty = VDV'$, where $D$ is the diagonal matrix of eigenvalues, and $V$ is the orthonormal matrix of eigenvectors. Let

$$T = D(D + \lambda I)^{-1}. \tag{15}$$

Then $\hat{f} = SY$, where $S = VTV'$. Mallows' $C_p$ in this context is $\|y - \hat{f}\|^2/n + (2/n)\operatorname{tr}(S)$. Notice the computation of the periodic Gaussian regularization in regression does not require equidistant design.

It is possible to leave the constant term in the regression function unpenalized, as is commonly done in practice with smoothing splines and Gaussian regularization. This is equivalent to having $\beta_0 = 0$ in (10), and the asymptotic results do not change. The penalized regression can be written as

$$\min_{f,b} \sum_{j=1}^n (y_j - (f(x_j) + b))^2 + \lambda J_0(f).$$

In this case the theory of reproducing kernel Hilbert spaces dictates that the solution can be expressed as $\hat{f} = G^\infty \hat{c} + \hat{b} e$, where $e = (1, \ldots, 1)'$. In the case of equidistant sample inputs, we can see that $e$ is an eigenvalue of $G^\infty$, since $G^\infty$ is periodic and even. In this case the computation is very similar to the case above with constants penalized: one simply replaces the diagonal element of $T$ in (15) corresponding to the eigenvalue $e$ by 1, and continues the computation with the new $T$.

**6. Simulations.** We conduct some simulations to study the finite sample properties of the periodic Gaussian regularization in regression. Consider the regression problem (2) with the following four functions on $[-\pi, \pi]$:

$$f_1(x) = \sin^2(x) \mathbb{1}_{(x \geq 0)},$$
$$f_2(x) = -x - \pi + 2(x + \pi/2) \mathbb{1}_{(x \geq -\pi/2)} + 2(-x + \pi/2) \mathbb{1}_{(x \geq \pi/2)},$$
$$f_3(x) = 1/(2 - \sin(x)),$$
$$f_4(x) = 2 + \sin(x) + 2\cos(x) + 3\sin^2(x) + 4\cos^3(x) + 5\sin^3(x).$$

The plots of the four functions are given in Figure 2. These are all $2\pi$-periodic functions. The first function has only the second order of smoothness. The second function has only the first order of smoothness. The third function is infinitely smooth. The fourth function is even smoother: its Fourier series only contains finitely many terms. In all of our simulations the sample size is taken to be 100. All simulations are done in Matlab.

First we study the effect of the joint tuning of $\lambda$ and $\omega$. We look at the regression problem (2) with the first regression function $f_1(x)$. In the first simulation we take the sample points to be equidistant in $(-\pi, \pi]$. The scatter plot is shown in Figure 3, top left panel. We use the periodic Gaussian regularization to do the estimation for $\omega = (k_1/5)^{(1/2)}$, $\lambda = \exp(-k_2/5)$, $k_1 = 1, \ldots, 100$, $k_2 = 1, \ldots, 100$. For each combination of $\omega$ and $\lambda$ we calculate the solution $\hat{f}_{\lambda,\omega}$ and the averaged squared error $(1/n)\sum_j [\hat{f}_{\lambda,\omega}(x_j) - f(x_j)]^2$. The bottom left panel of Figure 3 gives the corresponding contour plot of the averaged squared error. The $x$- and $y$-axes for the contour plot are $k_1$ and $k_2$, which are proportional to $\omega^2$ and $-\log \lambda$. Let the minimum of the averaged squared error be $a$. The levels in the contour plot are at $1.01a, 1.05a, 1.1a, 1.2a, 1.5a, 2a, 3a, 4a, 5a, 6a$. We used these levels to focus on the behavior of the averaged squared error around its minimum. It is clear that the contour levels are almost straight lines, indicating that the averaged squared errors are almost the same when $-\log \lambda$ varies linearly with $\omega^2$. This agrees with what is suggested by the asymptotic results, and suggests that in regression problems, as long as $\omega$ is fixed in a reasonable range, we can concentrate on the tuning of the smoothing parameter $\lambda$.

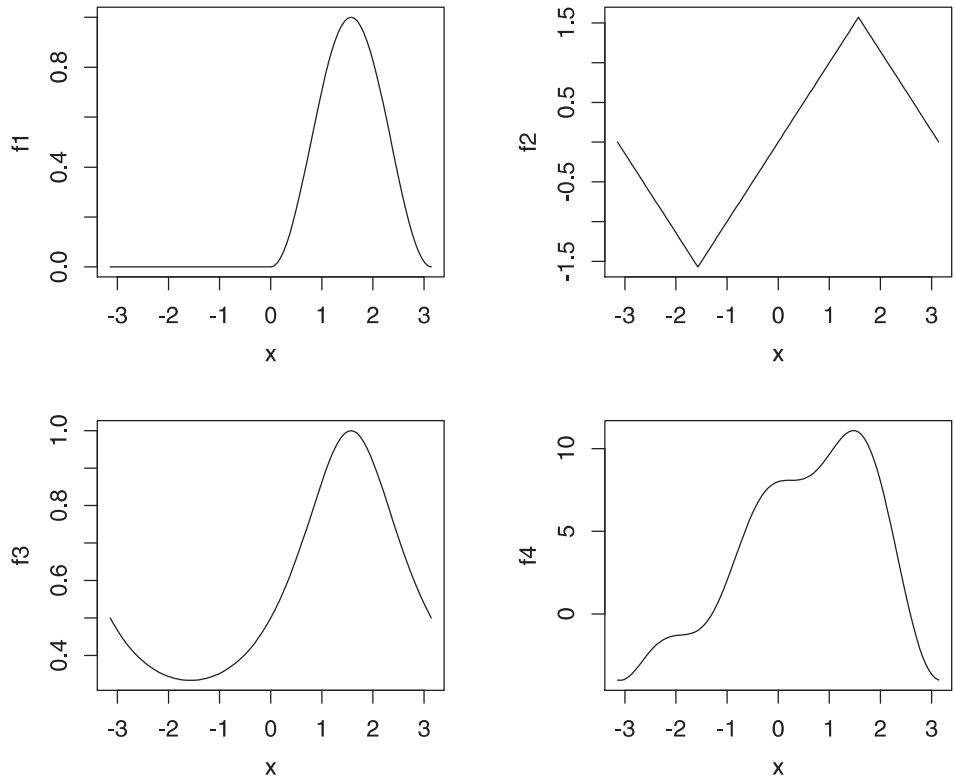

FIG. 2. *The regression functions used in the simulations. The first function has only the second order of smoothness. The second function has only the first order of smoothness. The third function is infinitely smooth. The fourth function has a Fourier series that only contains finitely many terms.*

Similar to any method of regularization, the periodic Gaussian regularization does not depend on the $x$'s being equidistant. The same phenomenon in the joint tuning of $\lambda$ and $\omega$ appears when the input $x$'s are not equidistant. We run the same simulation with nonequidistant $x$'s, and the corresponding scatter plot and the contour plot are given in the right panels of Figure 3. The nonequidistant $x$ values are generated by taking the fractional part of a normal variate with mean $1/4$ and standard deviation $1/4$, and then scaling the $[0, 1]$ interval to $[-\pi, \pi]$.

We run the same experiment with the other functions, $f_2$, $f_3$ and $f_4$, and the same observation about the joint tuning of $\lambda$ and $\omega$ is made in these experiments. This supports our strategy of a rough tuning for $\omega$ and a fine tuning over $\lambda$.

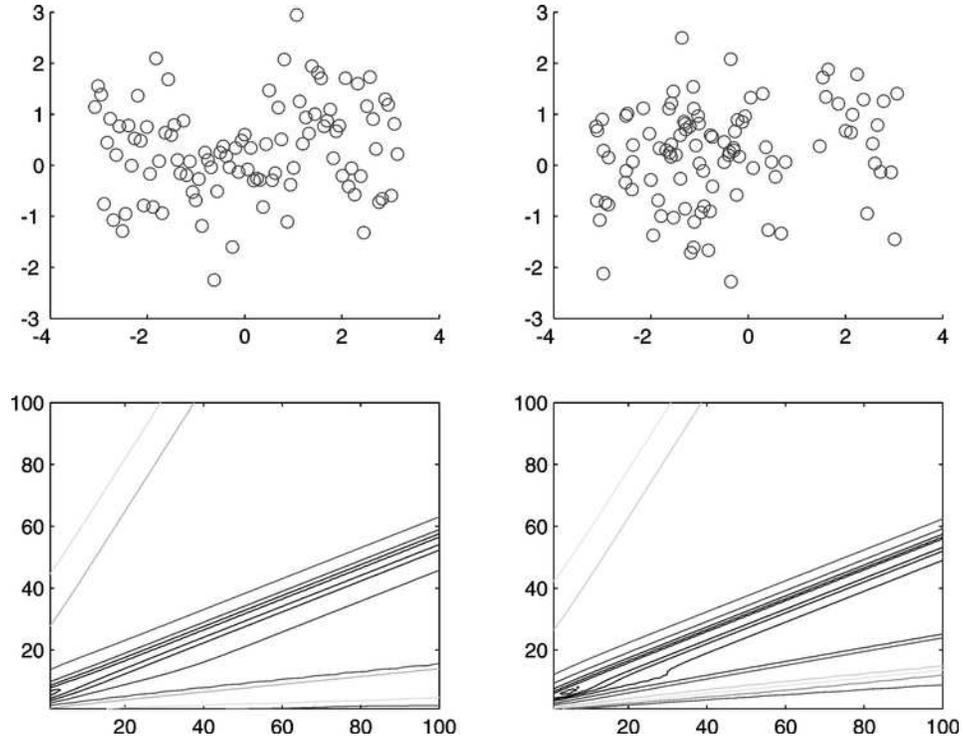

Fig. 3. *The top panels are the scatter plots of the data generated from the regression model* (2) *with the regression function* $f_1(x)$. Left: *equidistant case.* Right: *nonequidistant case. The bottom panels are the corresponding contour plots of the averaged squared errors of the periodic Gaussian regularization. The x- and y-axes for the contour plots are proportional to* $\omega^2$ *and* $-\log \lambda$, *respectively. We can see that in both cases the contour levels are very close to straight lines.*

TABLE 1
*Averaged squared error over 100 runs, for the periodic cubic smoothing spline, the periodic Gaussian regularization, and the periodic Gaussian regularization with constant left unpenalized, on four different functions of varying order of smoothness*

| Regression functions | Periodic cubic smoothing spline | Averaged squared error | |
|---|---|---|---|
| | | Periodic Gaussian regularization | |
| | | Constant penalized | Constant unpenalized |
| 1 | 0.0711 | 0.0675 | 0.0682 |
| 2 | 0.0541 | 0.0578 | 0.0582 |
| 3 | 0.0457 | 0.0462 | 0.0448 |
| 4 | 0.1136 | 0.0899 | 0.0899 |

Next we compare the periodic Gaussian regularization with the periodic cubic smoothing spline for regression on the circle on the four functions in Figure 2. The periodic cubic smoothing spline is the solution to

$$\sum_{j=1}^{n}(y_j - f(x_j))^2 + \lambda \int_{-\pi}^{\pi}[f''(t)]^2\,dt.$$

This penalty corresponds to the second-order Sobolev space, but leaves the linear functions unpenalized. For an introduction to the periodic cubic smoothing spline, see Wahba (1990) or Gu (2002).

We fix the $x$'s to be equidistant in $(-\pi, \pi)$ in our comparison. We generate $y$'s according to the regression model (2) with the four functions we consider. In both the periodic Gaussian regularization and the periodic cubic smoothing spline, the smoothing parameters are chosen according to Mallows' $C_p$. We search the minimal point of Mallows' $C_p$ over $\omega = 0.3k_1 - 0.1$, for $k_1 = 1, \ldots, 10$, and $\lambda = \exp(-0.4k_2 + 7)$, for $k_2 = 1, \ldots, 50$, for the periodic Gaussian regularization; and we search over $\lambda = \exp(-0.4k_2 + 7)$, for $k_2 = 1, \ldots, 50$, for the smoothing spline. We use the chosen smoothing parameter(s) to compute the solutions. For each generated dataset, we calculate the averaged squared error of the periodic Gaussian regularization and the periodic cubic smoothing spline.

We run the simulation 100 times. The averaged squared errors over the 100 runs are summarized in Table 1. For each regression function, a two-sided paired $t$-test is performed to compare the periodic Gaussian regularization and the periodic cubic smoothing spline based on the 100 runs. For the first function, the $p$-value is 0.49; for the second function, the $p$-value is 0.06, and it seems the smoothing spline may perform better; for the third function, the $p$-value is 0.9; for the fourth function, the $p$-value is very close to 0, and the periodic Gaussian regularization performed significantly better: we can see the averaged squared error of the periodic Gaussian regularization is 22% less than that of the periodic smoothing spline.

**7. Summary and discussion.** In this paper we study the method of regularization with the periodic Gaussian kernel. Asymptotically, the method adapts to unknown order of smoothness and is efficient compared with the minimax risk when the underlying function is reasonably smooth. The smoothing parameters in the periodic Gaussian regularization can be chosen adaptively without loss of asymptotic efficiency. Limited experiments in the finite sample case suggest that the performance of the periodic Gaussian regularization is comparable to that of the periodic cubic smoothing spline when the underlying regression function is reasonably smooth, and the periodic Gaussian regularization may have some advantage over the periodic cubic smoothing spline when the regression function is very smooth. This agrees with the asymptotic analysis, since it is well known that the cubic smoothing spline does not adapt to high order of smoothness.

The Gaussian reproducing kernel is commonly used in practice and has been successful in empirical studies. Our study on the periodic Gaussian reproducing kernel gives a partial explanation of the success of Gaussian reproducing kernel in practice, as we expect the Gaussian reproducing kernel to have similar properties to its periodic counterpart. When we apply the nonperiodic version of the Gaussian kernel to the examples in our simulation, the results are slightly inferior to the periodic version. This is to be expected, as the nonperiodic version does not take advantage of the fact that the functions in the simulation are periodic. However, the difference is not large. The averaged squared errors are 0.0736, 0.0679, 0.0559 and 0.1198.

The penalty functional $J_0$ in periodic Gaussian regularization corresponds to the norm of the infinite order Sobolev space $H_\omega^\infty$. It is also possible to consider the method of regularization with the penalty functional being the norm of the space $A_\alpha$ of analytic functions. This penalty cannot be written in terms of integrals of squared derivatives of integer order, but can be written in terms of derivatives of fractional order. In the Gaussian sequence model setting, the method of regularization with the analytic function space penalty is equivalent to the method of regularization (10) with $\beta_l = \exp(\alpha l)$. Similar asymptotic results as derived for the periodic Gaussian regularization can be derived for this alternative regularization: the method adapts to Sobolev space $H^m$ with unknown smoothness $m$. It is also possible to give an explicit expression for the reproducing kernel. In fact, the reproducing kernel is (14) with $\exp(-\omega^2 l^2/2)$ replaced by $\exp(-\alpha l)$. An equivalent form of this reproducing kernel is $E^\infty(s-t)$, with $E^\infty(r)$ defined as $E^\infty(r) \equiv \sum_{k=-\infty}^{\infty} E(r - 2k\pi)$ and $E(r) = \alpha/[\pi(r^2 + \alpha^2)]$ the Cauchy density function. This form follows from Lemma 3. Unlike the periodic Gaussian kernel case, the decay of $E(x)$ is slow, and it does not seem practical to use the form $E^\infty(s-t)$ for computation. On the other hand, it might be possible to calculate the reproducing kernel with the series in (14) with $\exp(-\alpha l)$.

## 8. Proofs.

PROOF OF THEOREM 1. The proof is an application of the theorem of Pinsker (1980). For completeness we state a form of the theorem given in Johnstone [(1998), Proposition 6.1 and Theorem 6.2]:

PINSKER'S THEOREM. *Consider the Gaussian sequence model* (9) *with the parameter space being the ellipsoid* $\Theta = \{\theta : \sum_l a_l^2 \theta_l^2 \leq Q\}$ *with* $a_l > 0$ *and* $a_l \to \infty$. *Then the minimax risk* $R(\Theta, n)$ *is asymptotically equivalent to the linear minimax risk* $R_L(\Theta, n)$, *which satisfies*

$$R_L(\Theta, n) = \frac{1}{n} \sum_l \left(1 - \frac{a_l}{\mu}\right)_+, \tag{16}$$

*where* $\mu = \mu(n, Q)$ *is determined by*

$$\frac{1}{n} \sum_l a_l (\mu - a_l)_+ = Q. \tag{17}$$

In our case we have $a_{2l} = a_{2l-1} = \exp(l^2 \omega^2 / 4)$, and (17) becomes

$$2 \sum_{l=1}^k \exp(l^2 \omega^2 / 4)\{\mu - \exp(l^2 \omega^2 / 4)\} = nQ,$$

with $k = k(\mu) = [2\omega^{-1}(\log \mu)^{1/2}]$, where $[\cdot]$ stands for the integer part. Notice that sums such as $\sum_{l=1}^k \exp(l^2 \omega^2 / 4)$ are dominated by the single leading term. Some calculations then give that $\log \mu \sim (1/2) \log(nQ)$. Therefore,

$$k = k(n) \sim 2^{1/2} \omega^{-1} (\log(nQ))^{1/2}.$$

Hence, it follows from Pinsker's theorem that

$$\begin{aligned} R(\Theta, n) &\sim R_L(\Theta, n) \\ &= \frac{1}{n} \sum_l \left(1 - \frac{a_l}{\mu}\right)_+ \\ &= \frac{2}{n} \sum_{l=1}^k \left\{1 - \frac{\exp(l^2 \omega^2 / 4)}{\mu}\right\} \\ &\sim \frac{2}{n} k(n) \sim 2^{3/2} n^{-1} \omega^{-1} (\log n)^{1/2}. \end{aligned}$$

This completes the proof of Theorem 1.
□

PROOF OF LEMMA 1. Solving the minimization problem (10), we get the method of regularization estimator $\hat{\theta}_l = (1 + \lambda\beta_l)^{-1} y_l$. As $\lambda$ goes to zero, we have

$$\sum_l \operatorname{var} \hat{\theta} = (1/n) \sum_l (1 + \lambda\beta_l)^{-2}$$

$$\sim (2/n) \sum_{l=1}^{\infty} (1 + \lambda e^{l^2 \omega^2/2})^{-2}$$

$$\sim (2/n) \int_0^\infty (1 + \lambda e^{x^2 \omega^2/2})^{-2} \, dx$$

$$= \sqrt{2} n^{-1} \omega^{-1} \int_{\log \lambda}^\infty (1 + e^y)^{-2} (y - \log \lambda)^{-1/2} \, dy$$

$$= \sqrt{2} n^{-1} \omega^{-1} \left[ \int_{\log \lambda}^0 (1 + e^y)^{-2} (y - \log \lambda)^{-1/2} \, dy \right.$$

$$\left. + \int_0^\infty (1 + e^y)^{-2} (y - \log \lambda)^{-1/2} \, dy \right].$$

For the second term in the bracket, we have

$$0 \le \int_0^\infty (1 + e^y)^{-2} (y - \log \lambda)^{-1/2} \, dy \le (-\log \lambda)^{-1/2} \int_0^\infty (1 + e^y)^{-2} \, dy.$$

Now let us look at the first term in the bracket. We have, on one hand,

$$\int_{\log \lambda}^0 (1 + e^y)^{-2} (y - \log \lambda)^{-1/2} \, dy \le \int_{\log \lambda}^0 (y - \log \lambda)^{-1/2} \, dy = 2(-\log \lambda)^{1/2};$$

on the other hand,

$$\int_{\log \lambda}^0 (1 + e^y)^{-2} (y - \log \lambda)^{-1/2} \, dy$$

$$\ge \int_{\log \lambda}^{-\log(-\log \lambda)} (1 + e^y)^{-2} (y - \log \lambda)^{-1/2} \, dy$$

$$\ge (1 + (-\log \lambda)^{-1})^{-2} \int_{\log \lambda}^{-\log(-\log \lambda)} (y - \log \lambda)^{-1/2} \, dy$$

$$\sim 2(-\log \lambda)^{1/2}.$$

Therefore, we have

$$\int_{\log \lambda}^0 (1 + e^y)^{-2} (y - \log \lambda)^{-1/2} \, dy \sim 2(-\log \lambda)^{1/2},$$

and the conclusion of the lemma follows. □

PROOF OF THEOREM 2. The periodic Gaussian regularization estimator is $\hat{\theta}_l = (1 + \lambda\beta_l)^{-1} y_l$. We have, for any $\theta \in H_\omega^\infty(Q)$,

$$\sum_{l=0}^\infty (E\hat{\theta}_l - \theta_l)^2 = \sum_{l=0}^\infty \lambda^2 \beta_l^2 (1 + \lambda\beta_l)^{-2} \theta_l^2$$

$$\leq 1/4\lambda \sum_{l=0}^\infty \beta_l \theta_l^2 = 1/4\lambda \sum_{l=0}^\infty \rho_l \theta_l^2 \leq 1/4\lambda Q.$$

Hence, from Lemma 1 we have, for any $\theta \in H_\omega^\infty(Q)$,

$$E \sum_l (\hat{\theta}_l - \theta_l)^2 = \sum_{l=0}^\infty (E\hat{\theta}_l - \theta_l)^2 + \sum_{l=0}^\infty \operatorname{var}\hat{\theta} \leq 2\sqrt{2}\omega^{-1} n^{-1} (-\log\lambda)^{1/2} + Q\lambda/4.$$

The last quantity is asymptotically equivalent to the asymptotic minimax risk $2\sqrt{2}\omega^{-1} n^{-1}(-\log\lambda)^{1/2}$ under (11). Therefore, under (11), the periodic Gaussian regularization estimator is asymptotically minimax. □

PROOF OF THEOREM 3. The estimator is $\hat{\theta}_l = (1 + \lambda\beta_l)^{-1} y_l$. From Lemma 1, we have

$$\sum_l \operatorname{var}\hat{\theta} = (1/n) \sum_l (1 + \lambda\beta_l)^{-2} \sim 2\sqrt{2}\omega^{-1} n^{-1} (-\log\lambda)^{1/2}.$$

On the other hand, we have

$$\sup_{\theta \in H^m(Q)} \sum_l (E\hat{\theta}_l - \theta_l)^2$$

$$= \sup_{\theta \in H^m(Q)} \sum_{l=0}^\infty (1 + \lambda^{-1}\beta_l^{-1})^{-2} \theta_l^2$$

$$= \sup_{\theta \in H^m(Q)} \sum_{l=0}^\infty (1 + \lambda^{-1}\beta_l^{-1})^{-2} \rho_l^{-1} (\rho_l \theta_l^2).$$

Here $\rho_{2l-1} = \rho_{2l} = 1 + l^{2m}$ are the coefficients in the definition (6) of the Sobolev ellipsoid $H^m(Q)$. Clearly, the maximum is achieved by putting all mass $Q$ at term $l$ that maximizes $(1 + \lambda^{-1}\beta_l^{-1})^{-2} \rho_l^{-1}$. That is, the maximum is

(18) $$Q\left[\max_l (1 + \lambda^{-1}\beta_l^{-1})^{-2} \rho_l^{-1}\right].$$

To evaluate (18), we first find the minimizer of

$$B_\lambda(x) = [1 + \lambda^{-1} \exp(-x^2\omega^2/2)]^2 (1 + x^{2m}) \qquad \text{over } x \geq 0.$$

Let $x_0(\lambda)$ be a global minimizer of $B_\lambda(x)$. It is easy to see that $x_0(\lambda) \neq \infty$, since $B_\lambda(\infty) = \infty$. Now let us first show that $x_0(\lambda) \to \infty$ as $\lambda \to 0$. We prove this with the elementary definition of limits. For any $M > 0$, we can find $\bar{x} > M$ such that $\exp[(\bar{x}^2 - M^2)\omega^2] > 1 + \bar{x}^{2m}$. Then $\lim_{\lambda \to 0} D(\lambda) > 1$, where

$$D(\lambda) = [\lambda + \exp(-M^2\omega^2/2)]^2 [\lambda + \exp(-\bar{x}^2\omega^2/2)]^{-2} (1 + \bar{x}^{2m})^{-1}.$$

Therefore, there exists $\delta > 0$, such that $D(\lambda) > 1$ for any $\lambda < \delta$. On the other hand, for any $x \leq M$, we have $B_\lambda(x)/B_\lambda(\bar{x}) \geq D(\lambda)$. Therefore, for any $\lambda < \delta$, we have $\inf_{x \leq M} B_\lambda(x) > B_\lambda(\bar{x})$, therefore, $x_0(\lambda) > M$. This shows that $x_0(\lambda) \to \infty$ as $\lambda \to 0$.

Since $x_0(\lambda) \neq \infty$, we have $B'_\lambda(x_0) = 0$. That is,

$$\tag{19} m^{-1}\omega^2(x_0^2 + x_0^{-(2m-2)}) = 1 + \lambda \exp(x_0^2\omega^2/2).$$

Since $x_0(\lambda) \to \infty$ as $\lambda \to 0$, we have

$$\tag{20} m^{-1}\omega^2 x_0^2 \sim \lambda \exp(x_0^2\omega^2/2),$$

$$\tag{21} x_0^2\omega^2/2 \sim (-\log \lambda).$$

Therefore, by (19) and (20) we have

$$B_\lambda(x_0) = [1 + (m^{-1}\omega^2(x_0^2 + x_0^{-(2m-2)}) - 1)^{-1}]^2 (1 + x_0^{2m}) \sim x_0^{2m}.$$

From this and (21), we see that

$$Q\left[\max_l (1 + \lambda^{-1}\beta_l^{-1})^{-2} \rho_l^{-1}\right] \sim Q x_0^{-2m} \sim Q 2^{-m} \omega^{2m} (-\log \lambda)^{-m}.$$

Therefore,

$$\tag{22} \max_{\theta \in H^m(Q)} \sum_l E(\hat{\theta}_l - \theta_l)^2$$
$$\sim Q 2^{-m} \omega^{2m} (-\log \lambda)^{-m} + 2\sqrt{2} \omega^{-1} n^{-1} (-\log \lambda)^{1/2}.$$

The conclusion of the theorem then comes from simple calculations. □

PROOF OF THEOREM 5. By (13), we have

$$\sup_{\theta \in H^m(Q)} \sum_l E(\hat{\theta}_l^* - \theta_l)^2$$
$$\leq (1 + O(B^{-1})) \sup_{\theta \in H^m(Q)} \min_{\lambda, \omega_s} \left\{ \sum_l E(\hat{\theta}_l - \theta_l)^2 \right\} + n^{-1} O(B)$$
$$\leq (1 + O(B^{-1})) \min_{\lambda, \omega_s} \sup_{\theta \in H^m(Q)} \left\{ \sum_l E(\hat{\theta}_l - \theta_l)^2 \right\} + n^{-1} O(B).$$

Similar inequalities hold for $H^\infty_{\omega_s}(Q)$ and $A_\alpha(Q)$. Now take $B = (\log n)^{1/3}$, and the conclusion of the theorem follows from Theorems 1–4 and the fact that $\omega_s \in \{\omega_1, \ldots, \omega_S\}$ has finitely many possibilities. □

Department of Statistics  
University of Wisconsin  
1210 West Dayton street  
Madison, Wisconsin 53706  
USA  
e-mail: yilin@cs.wisc.edu

Department of Statistics  
University of Pennsylvania  
3730 Walnut Street  
Philadelphia, Pennsylvania 19104-6340  
USA  
e-mail: lbrown@wharton.upenn.edu